\documentclass[12pt]{article}
\usepackage{amssymb,latexsym}
\usepackage{amsmath,bbm,amsthm}
\usepackage[active]{srcltx}
\def\finproof{\hbox{\vrule width1.0ex height1.ex}\vspace{2mm}}
\begin{document}

%%%%%%%%%%%%%%%%%%%%   without "ignoring" term  %%%%%%%%%%%%%%%%%%%%

\begin{center}
{\large \sc A condition on delay for differential equations with
discrete state-dependent
delay}%: state-dependent ignoring condition }%

\bigskip

{\sc A.V.Rezounenko}

\smallskip

Department of Mechanics and Mathematics \\ V.N.Karazin Kharkiv
National University \\ 4, Svobody Sqr., Kharkiv, 61077, Ukraine \\
{\it E-mail}: rezounenko@univer.kharkov.ua

\end{center}

\noindent {\bf Abstract.} Parabolic differential equations with
discrete state-dependent delay are studied. The approach, based on
an additional condition on the delay function introduced in [A.V.
Rezounenko, Differential equations with discrete state-dependent
delay: uniqueness and well-posedness in the space of continuous
functions, Nonlinear Analysis: Theory, Methods and Applications, 70
(11) (2009), 3978-3986] is developed. We propose and study a {\it
state-dependent} analogue of the condition which is sufficient for
the well-posedness of the corresponding initial value problem on the
whole space of continuous functions $C$. The dynamical system is
constructed in $C$ and the existence of a compact global attractor
is proved.\\ %

\smallskip

\noindent {\it AMS subject classification}: 35R10 35B41 35K57. \\ %

\smallskip

\noindent {\it Keywords}: Partial functional differential equation,
State-dependent delay, Well-posedness, Global attractor.

\section{Introduction }\label{sec1}
%\marginpar{\tiny sdd12}

%\marginpar{\tiny 2010.11.17}
Delay differential equations is one of
the oldest branches of the theory of infinite dimensional dynamical
systems - theory which describes qualitative properties of systems,
changing in time.

We refer to the classical monographs on the theory of ordinary
(O.D.E.) delay equations
\cite{Hale,Hale_book,Walther_book,Azbelev,Mishkis}. The theory of
partial (P.D.E.) delay equations is essentially less developed
since such equations are infinite-dimensional in both time (as
delay equations) and space  (as P.D.E.s) variables, which makes
the analysis more difficult. We refer to some works which are
close to the present research
\cite{%travis_webb,
Chueshov-JSM-1992,Cras-1995,NA-1998,Rezounenko-2003}
and to the monograph \cite{Wu_book}.

 A new class of equations with delays has recently attracted
attention of many researchers. These equations have a delay term
that may depend on the state of the system, i.e. the delay is {\tt
state-dependent} (SDD). Due to this type of delays such equations
are inherently nonlinear and their study has begun in the case of
ordinary differential equations
\cite{Nussbaum-Mallet-1992,Nussbaum-Mallet-1996,MalletParet,Krisztin-Arino-2001,
Walther2002,Walther_JDE-2003,Krisztin-2003} (for more details see
also a recent survey \cite{Hartung-Krisztin-Walther-Wu-2006},
articles \cite{Walther_JDDE-2007,Walther-JDDE-2010} and references
therein).

Investigations of these equations essentially differ from the ones
of equations with constant or time-dependent delays. The underlying
main mathematical
difficulty of the theory %of equations %PDEs with (discrete) SDDs
lies in the fact that delay terms %the functions describing
with discrete state-dependent delays
%(in contrast to constant or time-dependent delays)
are not Lipschitz continuous on the space of continuous functions -
the main space, on which the classical theory of equations with
delays is developed  (see \cite{Winston-1970} for an explicit
 example of the non-uniqueness and
 \cite{Hartung-Krisztin-Walther-Wu-2006} for more details). It is a common point of
view~\cite{Hartung-Krisztin-Walther-Wu-2006} that the corresponding
\textit{initial value problem} (IVP) is not generally well-posed in
the sense of J.~Hadamard \cite{Hadamard-1902,Hadamard-1932} in the
space of continuous functions ($C$). This leads to the search of
(particular) classes of equations which may be well-posed in the
space of continuous functions ($C$).

Results for partial differential equations with SDD  %in this direction
have been %were
obtained only recently in
 \cite{Rezounenko-Wu-2006}%
 (case of distributed delays, weak solutions), %{}
\cite{Hernandez-2006} (mild solutions, unbounded discrete delay),
and \cite{Rezounenko_JMAA-2007} (weak solutions, bounded discrete
and distributed delays).

The main goal of the present paper is to develop an alternative
approach, based
on an additional condition (see (H) below)   %the so-called "ignoring condition"
introduced in
\cite{Rezounenko_NA-2009}. We propose and study a {\it
state-dependent} analogue of the condition %"ignoring condition"
which is sufficient for the
well-posedness of the corresponding initial value problem in the
space C. The presented approach includes the possibility when the
state-dependent delay function does not satisfy the %ignoring
condition on a subset of the phase space $C$, but the IVP still be
well-posed in the whole space $C$. This is our second goal which is
to connect the approach developed for ODEs (a restriction to a
subset of Lipschitz continuous functions) and the
approach~\cite{Rezounenko_NA-2009} of a different nature.

Discussing the meaning of the main %"ignoring"
assumptions (H) and
$(\widehat{H})$ (see below) for
 %many
applied problems, we hope that these assumptions are the natural
mathematical expression of the fact that %for many applied problems,
%%%%%%%%%%%%%%%%%%%%%%%%%%%%%%%%%%%%%
many differential equations encountered in modeling real world
phenomena
%%%%%%%%%%%%%%%%%%%%%%%%%%%%%%%%%%%%%
%the delay time is
%the models
have a parameter (time $\eta_{ign}>0$ or $\Theta^\ell > 0$) which is
necessary to take into considerations the time changes in the
system. The changes not always can be taken into considerations
immediately. To this end, the existence of $\eta_{ign}>0$ or
$\Theta^\ell > 0$ (no matter how small the values of $\eta_{ign}>0$
or $\Theta^\ell > 0$ are!) makes the corresponding initial value
problem well-posed.

Having the well-posedness proved, we study the long-time
asymptotic behavior of the correspond dynamical system and prove
the existence of a compact global attractor.

\section{Formulation of the model with state-dependent discrete delay}\label{sec2}
%\marginpar{\tiny new}

Let us consider the following  %\marginpar{\tiny 2010.11.11} %
%non-local
parabolic partial differential equation with %{\it state-dependent discrete %\underline{discrete}
delay
%}
%\begin{equation}\label{sdd3-1}
%\begin{array}{lll}
%&\,\,\,\, \frac{\partial }{\partial t}u(t,x)+Au(t,x)+du(t,x)\\
%&= %\varepsilon
%%\int^0_{-r} \left\{
%\int_\Omega b(u(t-\eta (u(t),u_t), y)) f(x-y) dy   \equiv \big(
%F(u_t) \big)(x),\quad x\in \Omega ,\end{array}
%\end{equation}

\begin{equation}\label{sdd8-1}
\frac{\partial }{\partial t}u(t,x)+Au(t,x)+du(t,x)
%\int^0_{-r} \left\{
=
%\int_\Omega b(u(t-\eta (u_t), y)) f(x-y) dy \equiv
\big( F(u_t)
\big)(x),\quad x\in \Omega ,%\end{array}
\end{equation}
 where $A$ is a densely-defined self-adjoint positive linear operator
 with domain $D(A)\subset L^2(\Omega )$ and with compact
  resolvent, so $A: D(A)\to L^2(\Omega )$ generates an analytic semigroup,
  $\Omega $ is a smooth bounded domain in $R^{n_0}$, $d$ is a non-negative constant.
As usually for delay equations, we denote by $u_t$ the function of
$\theta\in [-r,0]$ by the formula $u_t\equiv u_t(\theta)\equiv
u(t+\theta).$ We denote for short $C\equiv C([-r,0];L^2(\Omega)).$
The norms in $L^2(\Omega)$ and $C$ are denoted by $||\cdot ||$ and
$||\cdot ||_C$ respectively.

\smallskip

The (nonlinear) delay term $F : C([-r,0];L^2(\Omega)) \to
L^2(\Omega)$ has the form
\begin{equation}\label{sdd12-2-1}
  F(\varphi) = B(\varphi (-\eta (\varphi))),
\end{equation}
where (nonlinear) mapping $B : L^2(\Omega) \to L^2(\Omega)$ is
Lipschitz continuous
\begin{equation}\label{sdd12-2-2}
||B (v^1) - B(v^2)|| \le L_B ||\, v^1 - v^2||,\quad \forall v^1,v^2
\in L^2(\Omega).
\end{equation}

%$f: \Omega -\Omega \to R$   is a bounded
%function to   be specified later,
%$b:R\to R$ is a locally Lipschitz %bounded
%map, %  ($|b(w)|\le C_b$ with $C_b\ge 0),$
%
%

\smallskip

 The function
  $\eta (\cdot): C([-r,0];L^2(\Omega)) \to [0,r]\subset R_{+}$ represents the
{\tt state-dependent discrete delay}.
%We denote for short $H\equiv L^2(\Omega)\times L^2(-r,0;L^2(\Omega)).$
It is important to notice that $F$ is {\it nonlinear} even in the
case of linear~$B$.

We consider equation (\ref{sdd8-1}) with the following initial
condition
\begin{equation}\label{sdd8-ic}
u|_{[-r,0]}=\varphi \in C\equiv C([-r,0];L^2(\Omega)).
\end{equation}

\medskip

{\bf Remark~1.} {\it The results presented in this paper could be
easily extended to the case of nonlinearity $F$ %(see (\ref{sdd12-2-1}))
 of the form $F(\varphi) = \sum_k B^k(\varphi
(-\eta^k (\varphi)))$ as well as to O.D.E.s, for example, of the
following form~\cite{Rezounenko-metod-2010}
%\medskip
\begin{equation}\label{sdd8-5-12}
\dot u (t) + A u(t)+ d\cdot u(t) = b(u(t-\eta (u_t))), \quad
u(\cdot)\in R^n, \, d\ge 0.
\end{equation}
In the last case one simply needs to substitute $L^2(\Omega )$ by
$R^n$ and use $C\equiv C([-r,0];R^n)$ instead of  $
C([-r,0];L^2(\Omega )).$
%Assumptions on the delay function $\eta$ are the same.
The function $b : R^n\to R^n$ is locally Lipschitz continuous and
satisfies $||b(w)||_{R^n} \le C_1 ||w||_{R^n} +C_b$
with $C_1, C_b \ge 0;$ $A$ is a matrix. %
}%

\medskip

{\bf Remark~2.} {\it As an example we could consider {\tt nonlocal}
delay term $F$ (see (\ref{sdd12-2-1})) with the following mapping
$$B(v)(x)\equiv \int_\Omega b(v(y)) f(x-y) dy, \quad x\in \Omega,$$
where $f: \Omega -\Omega \to R$ is a bounded and measurable function
($|f(z)|\le M_f, \forall z\in \Omega -\Omega$) and $b:R\to R$ is
a (locally) Lipschitz mapping, %bounded
 %  ($|b(w)|\le C_b$ with $C_b\ge 0),$
satisfying $|b(w)|\le C_1|w|+C_b$ with $C_i\ge 0.$ In this case
equation (\ref{sdd8-1}) has the form
$$\frac{\partial }{\partial t}u(t,x)+Au(t,x)+du(t,x) = \int_\Omega b(u(t-\eta (u_t), y)) f(x-y)
dy,\quad x\in \Omega.
$$
One can easily check that $B$ satisfies (\ref{sdd12-2-2}) with
$L_B\equiv L_b M_f |\Omega|$, where $L_b$ is the Lipschitz constant
of $b$,  and $|\Omega| \equiv \int_\Omega 1 \, dx$.

Another example is a ({\tt local}) delay term $F$ (see
(\ref{sdd12-2-1}))
with %the following
 $B(v)(x)\equiv  b(v(x)), x\in \Omega.$ Equation (\ref{sdd8-1}) has the
 form
$$\frac{\partial }{\partial t}u(t,x)+Au(t,x)+du(t,x) =  b(u(t-\eta (u_t), x)), \quad x\in
\Omega.
$$
An easy calculation show that (\ref{sdd12-2-2}) is satisfied with
$L_B\equiv L_b$.
}%

\medskip

The methods used in our work can be applied to other types of
nonlinear and delay P.D.E.s (as well as O.D.E.s). We choose a
particular form of nonlinear delay terms $F$ for simplicity and to
illustrate our approach on the diffusive Nicholson's blowflies
equation (see the end of the article for more details).

%{\bf Remark.}

\section{The existence of mild solutions}\label{sec3}

In our study we use the standard

\smallskip

{\bf Definition~1}. %\marginpar{?} %
 {\it A function $u\in C([-r,T]; L^2(\Omega))$ is called a {\tt mild solution}
 on $[-r,T]$ of the initial value problem (\ref{sdd8-1}), (\ref{sdd8-ic}) if it satisfies
 (\ref{sdd8-ic}) and
 \begin{equation}\label{sdd8-3-1}
u(t)=e^{-A t}\varphi(0) + \int^{t}_0 e^{- A (t-s)} \left\{ F(u_s)
- d \cdot u(s)\right\}\, ds, \quad t\in [0,T].
 \end{equation}
}%

\medskip

{\bf Proposition~1}\cite{Rezounenko_NA-2009}. {\it Assume
%the function $b : R\to R$ is a locally Lipschitz %bounded
%map, satisfying $|b(w)|\le C_1|w|+C_b$ with $C_i\ge 0,$ and
%  $f: \Omega -\Omega \to R$   is a bounded and measurable function.
the mapping $B$ is Lipschitz continuous (see (\ref{sdd12-2-2})) and
 delay  function
  $\eta (\cdot): C([-r,0];L^2(\Omega)) \to [0,r]\subset R_{+}$ is
  continuous.

  Then for any
  initial function $\varphi\in C,$ initial value problem (\ref{sdd8-1}), (\ref{sdd8-ic})
  has a global mild solution which satisfies $u\in C([-r,+\infty); L^2(\Omega))$.
  }%

\medskip

%The existence of solutions in $C.$

The existence of a mild solution is a consequence of the continuity
of $F
%$F(\varphi)\equiv\int_\Omega b(\varphi (-\eta
%(\varphi), y)) f(\cdot -y) dy
% %%%b(\varphi(-\eta (\varphi)))
: C\to L^2(\Omega)$ (see (\ref{sdd8-1})) which gives the possibility
to use the standard method based on Schauder fixed point theorem
(see e.g. \cite[theorem 2.1, p.46]{Wu_book}). The solution is also
global (is defined for all $t\ge -r$) since (\ref{sdd12-2-2})
implies $||F(\varphi)||\le L_B ||\varphi||_C + ||B(0)||$ and one can
apply, for example, \cite[theorem 2.3, p.49]{Wu_book}.

%
%see above and th 1 [Kantorovich-Akilov] (p.424, Part XI, section 3)
% on the "measurable".
%

\medskip

{\bf Remark~3.} {\it It is important to notice that even in the case
of ordinary differential equations (even scalar) the mapping of the
form $\widetilde{F} (\varphi)=\widetilde f(\varphi(-r(\varphi)))\,
:\, C([-r_0,0];R) \to R$ has a very unpleasant property. The authors
in \cite[p.3]{Louihi-Hbid-Arino-JDE-2002} write "Notice that the
functional $\widetilde{F}$ is defined on $C([-r_0,0];R),$ but it is
clear that it is neither differentiable nor locally Lipschitz
continuous, whatever the smoothness of $\widetilde f$ and $r.$" As a
consequence, the Cauchy problem associated with equations with such
a nonlinearity "...is \texttt{not} well-posed in the space of
continuous functions, due to the non-uniqueness of solutions
whatever the regularity of the functions $\widetilde f$ and $r$"
\cite[p.2]{Louihi-Hbid-Arino-JDE-2002}. See also a detailed
discussion in \cite{Hartung-Krisztin-Walther-Wu-2006}.%
}%

\medskip

{\bf Remark~4.} {\it For a study of solutions to %PDEs
equations with a state-dependent delay in the space $C([-r,0];E) $
with $E$ not necessarily finite-dimensional Banach space see e.g.
\cite{Arino-Sanchez-DCDS-2005}
% could be proved for a wider class of equations (see e.g. \cite{}).

}%

\medskip

%Our main goal in
In this work we concentrate on conditions for the IVP
(\ref{sdd8-1}), (\ref{sdd8-ic}) to be well-posed.

\section{Main results: uniqueness, well-posedness and asymptotic bahavior}\label{sec4}

As in the previous section, we assume that $\eta : C \to [0,r]$ is
continuous and $B$ is Lipschitz.
%  $f: \Omega -\Omega \to R$   is a bounded and
%measurable function ($|f(\cdot )|\le M_f$).
Unlike to the existence of solutions, the uniqueness is essentially
more delicate question in the presence of discrete state-dependent
delay (see a classical example of the non-uniqueness in
\cite{Winston-1970}).

Let us remind an important additional assumption  on the delay
function
$\eta$, %the so-called {\bf "ignoring condition"}
as it was
introduced in \cite{Rezounenko_NA-2009}:

\begin{itemize}
 \item $\exists \eta_{ign}>0$ such that $\eta$ "ignores" values of
 $\varphi(\theta)$ for $\theta\in (-\eta_{ign},0]$ i.e. %for any
 $$\hskip-10mm \exists\, \eta_{ign}>0 : %\quad
 \forall\varphi^1, \varphi^2\in C :
 \forall\theta\in
 [-r,-\eta_{ign}],\,\Rightarrow \varphi^1(\theta)= \varphi^2(\theta)\quad
   \Longrightarrow \quad
 \eta (\varphi^1)=\eta (\varphi^2). \eqno(H) $$
 %%%%\eqno(H2)$$
\end{itemize}
\medskip

For examples of delay functions satisfying (H) and the proof of the
{\it uniqueness} of mild solutions (given by Proposition~1) as well
as the well-posedness of the IVP (\ref{sdd8-1}), (\ref{sdd8-ic}) see
\cite{Rezounenko_NA-2009}.

\medskip

{\bf Remark~5.} {\it It is important  to notice that, discussing the
%ignoring
condition (H) and its dependence on the %ignoring
value $\eta_{ign}$,
we see that  in the case $\eta_{ign}>r,$ one has that the delay
function $\eta$ ignores {\tt all} values of $\varphi(\theta),
\forall \theta\in [-r,0],$ so $\eta (\varphi)\equiv const,
\forall\varphi\in C$ i.e. equation (\ref{sdd8-1}) becomes an
equation with {\tt constant (!)} delay. On the other hand, the
analogue of assumption (H) with $\eta_{ign}=0,$ is trivial since
$\varphi^1(\theta)=\varphi^2(\theta)$ for all $\theta\in [-r,0]$
means $\varphi^1=\varphi^2$ in $C,$ so $\eta (\varphi^1)=\eta
(\varphi^2).$ }

\medskip

{\bf Remark~6.} {\it It is worth mentioning that the classical case
of {\tt constant} delay (see the previous remark) and the
corresponding theory forms the basement for the discussed approach,
%(ignoring condition approach),
but could be mixed with the approach of non-vanishing delays. In our
case the delay $\eta$ do may vanish (we do {\bf not} assume the
existence of $r_0>0$ such that $\eta(\varphi)\ge r_0, \forall
\varphi$).
}%

 \medskip

In the above %{\bf "ignoring condition"}
condition (H) the %ignored
semi-interval $(-\eta_{ign},0]$ is fixed (we remind that the value
%{\bf "ignoring value"}
 $\eta_{ign}$ could be arbitrary small).

Our goal is to extend the approach based on the condition (H) %"ignoring condition"
to a more wide class of state-dependent delay functions where the
value %{\bf "ignoring value"}
$\eta_{ign}$ is not a constant any more, but a function of the
state. Moreover, as an easy additional extension, we also allow the
upper bound of the delayed segment to be state-dependent. More
precisely, we consider two functions $\Theta^u,
\Theta^\ell : C\to [0,r]$ (upper and low %ignoring
functions),
satisfying
$$\forall \varphi\in C \quad \Rightarrow \quad 0\le \Theta^\ell  (\varphi) \le
\Theta^u
(\varphi)\le r.
$$
Now we are ready to introduce~\cite{Rezounenko-metod-2010} the following {\it state-dependent} %ignoring
condition for the state-dependent delay function $\eta:C\to [0,r]$
(c.f. (H)):

\begin{itemize}
 \item  $\eta$ "ignores" values of
 $\varphi(\theta)$ for $\theta\not\in [-\Theta^u  (\varphi),
-\Theta^\ell (\varphi)]$ i.e. %for any
 $$\hskip-10mm
 \forall\, \psi \in C  \mbox{ such that }\,
 \forall\theta\in
 [-\Theta^u  (\varphi),
-\Theta^\ell (\varphi)]\,\Rightarrow \psi(\theta)=
\varphi(\theta)\quad
   \Longrightarrow \quad
 \eta (\psi)=\eta (\varphi). \eqno(\widehat{H}) $$
 %%%%\eqno(H2)$$
\end{itemize}

%\begin{itemize}
% \item  $\eta$ "ignores" values of
% $\varphi(\theta)$ for $\theta\not\in [-\eta^u  (\varphi),
%-\eta^\ell (\varphi)]$ i.e. %for any
% $$\hskip-10mm %\exists\, \eta_{ign}>0 : %\quad
% \forall\varphi^1, \varphi^2\in C : %\mbox{such that}
% \forall\theta\in
% [-\eta^u  (\varphi),
%-\eta^\ell (\varphi)]\,\Rightarrow \varphi^1(\theta)=
%\varphi^2(\theta)\quad
%   \Longrightarrow \quad
% \eta (\varphi^1)=\eta (\varphi^2). \eqno(\widehat{H}) $$
%\end{itemize}
%

\medskip

The above condition means that state-dependent delay function $\eta$
"ignores" all values of its argument $\varphi$ outside of
$[-\Theta^u (\varphi), -\Theta^\ell (\varphi)]\subset [-r,0]$ and
this delayed segment $[-\Theta^u  (\varphi), -\Theta^\ell
(\varphi)]$ is {\tt state-dependent}. We could illustrate this
property on the picture.
%\marginpar{\tiny PICTURE !}

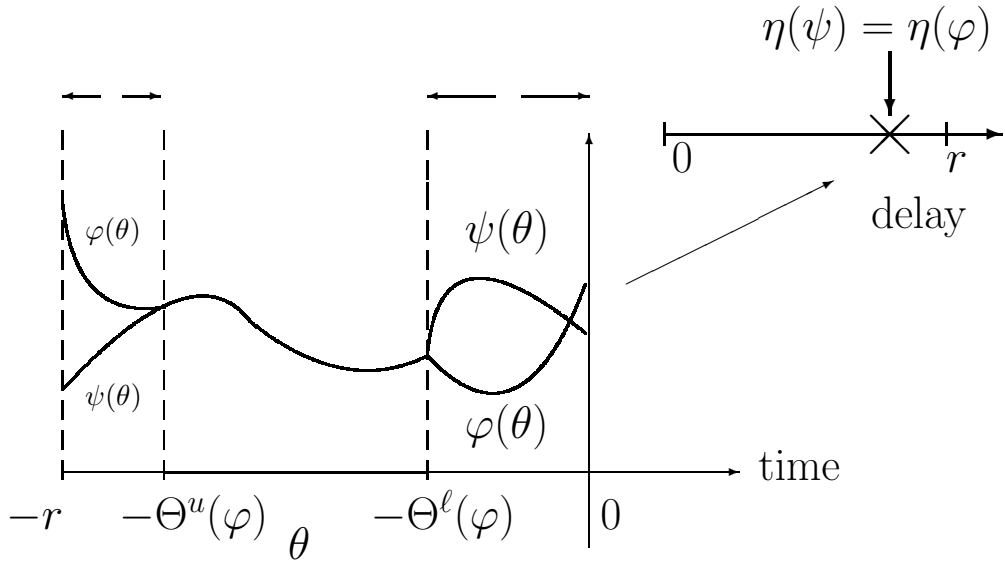
\begin{Large}
\begin{center}
\hskip70mm
\begin{picture}(400,200)%%%%%%%%begin-ignoring condition-illustration%%%%%%%%%%%%%%%%%%%%%%%%%%
\unitlength=0.5mm%mashtab picture!!!
 \put(0,20){\vector(1,0){180}}
\put(140,0){\vector(0,1){110}} \put(185,18){\text{time}}

\put(143,5){$0$}\put(140,18){\line(0,1){5}}
%\put(67,0){\Large $s-\varepsilon_k$}\put(70,18){\line(0,1){5}}
\put(97,90){\line(0,1){7}}

%\put(20,23){\Large $s-\varepsilon_n$}

%\put(70,18){\line(0,1){5}}
%\put(145,75){$1/\varepsilon_n$}\put(138,80){\line(1,0){5}}
%\put(145,95){$1/\varepsilon_k$}\put(138,97){\line(1,0){5}}

\multiput(0,18)(0,8){12}{\line(0,1){5}}%right
\put(-15,5){$-r$}\put(0,18){\line(0,1){5}}
%\multiput(70,18)(0,18){5}{\line(0,1){7}}%center
%\multiput(57,18)(0,8){8}{\line(0,1){5}}%left
\multiput(97,18)(0,8){12}{\line(0,1){5}}%right
\put(82,5){$-\Theta^\ell (\varphi)$}
\multiput(27,18)(0,8){12}{\line(0,1){5}}%right
\put(15,5){$-\Theta^u (\varphi)$}
 \put(60,-3){$\theta$}
%\put(-20,60){$\widetilde\xi^n(\theta,s)$}
%\put(-5,100){$\widetilde\xi^k(\theta,s)$}

%%%%%%%%%%%%%%%%common part of \varphi^1 and \varphi^2 %%%%%%%%%%%%
\put(0,40){\begin{picture}(60,60)%%%%%%%%%%%%%
%\textcolor[rgb]{0.00,1.00,0.00}{ %%%green color%%%%
 \thicklines
\qbezier(0,2)(35,40)(50,20)\qbezier(50,20)(75,0)(97,11)
%}%

\end{picture}}%%%%%%%%%%%%%%%%%%%%%%%%%%%%%%%%%
%%%%%%%%%%%%%%%%a part1 of \varphi^1 %%%%%%%%%%%%
\put(0,40){\begin{picture}(60,60)%%%%%%%%%%%%%
{\thicklines \qbezier(0,52)(3,20)(27,24)}%
\put(6,42){\small $\varphi(\theta)$}
\end{picture}}%%%%%%%%%%%%%%%%%%%%%%%%%%%%%%%%%
%%%%%%%%%%%%%%%%a part2 of \varphi^1 %%%%%%%%%%%%
\put(97,49){\begin{picture}(60,60)%%%%%%%%%%%%%
{\thicklines \qbezier(0,2)(3,40)(42,8)}%
\put(10,32){$\psi(\theta)$}
\end{picture}}%%%%%%%%%%%%%%%%%%%%%%%%%%%%%%%%%
%%%%%%%%%%%%%%%%a part1 of \varphi^2 %%%%%%%%%%%%
\put(6,38){\small $\psi(\theta)$}
%%%%%%%%%%%%%%%%a part2 of \varphi^2 %%%%%%%%%%%%
\put(97,49){\begin{picture}(60,60)%%%%%%%%%%%%%
{\thicklines \qbezier(0,2)(25,-25)(42,21)}%\qbezier(0,2)(20,10)(42,95)
\put(10,-20){$\varphi(\theta)$}
\end{picture}}%%%%%%%%%%%%%%%%%%%%%%%%%%%%%%%%%

%%%%%%%%%%%%%%%%\eta( \varphi^1)=\eta(\varphi^2) %%%%%%%%%%%%
\put(160,90){\begin{picture}(60,60)%%%%%%%%%%%%%
\thicklines \put(0,20){\vector(1,0){90}}

\put(2,10){$0$}\put(0,17){\line(0,1){6}}
\put(75,10){$r$}\put(75,17){\line(0,1){6}}%
\put(55,15){\line(1,1){10}} \put(55,25){\line(1,-1){10}}
\put(26,45){$\eta( \psi)=\eta(\varphi)$}
\put(60,42){\vector(0,-1){17}} \put(55,-5){\text{delay}}
\end{picture}}%%%%%%%%%%%%%%%%%%%%%%%%%%%%%%%%%

%%%%%red line under \varphi^1 and \varphi^2 %%%%%%%%%%%%
{\thicklines
{\put(27,20){\line(5,0){70}}}%
}%%%%%%%%%%%%%%%%%%%%%%%%%%%%%%%%%%%%%%%%%%%%%%%%%%%%%%%%%

\put(-10,130){{\text{ %ignored
}}} %%
{ %left ignored set
\put(16,120){\vector(1,0){10}}%%
\put(10,120){\vector(-1,0){10}}%%
 }%%

\put(95,130){{\text{ %ignored
}}} %%
{ %right ignored set
\put(122,120){\vector(1,0){18}}%%
\put(115,120){\vector(-1,0){18}}%%
 }%%
%\put(95,-10){\text{ ignored}}

%%%%%%%%%%%%%%%%%%%%%%%%%%
\put(150,70){\vector(2,1){55}}%
%%%%%%%%%%%%%%%%%%%%%%%%arrow{}

\end{picture}%%%%%%%%4%%%%%%%%%%%%%%%%%%%%%%%%%%
%\end{right}%
\end{center}

\end{Large}

\bigskip

\medskip

{\bf Remark~7.} {\it One could see that (H) is a particular case of
$(\widehat{H})$ with $\Theta^\ell (\varphi)\equiv \eta_{ign}$ and
$\Theta^u  (\varphi)\equiv r, \forall\varphi\in C.$}

\medskip

{\bf Examples.} {\it It is easy to present many examples of (delay)
functions $\eta$, which satisfy assumption $(\widehat{H})$. The
simplest one is

\begin{equation}\label{sdd12-4-4}
%$$
\eta(\varphi)= p_1\left( \varphi(-\chi(\varphi(-r))\right) \hbox{
with } p_1: L^2(\Omega)\to [0,r]
%$$
\end{equation}
and given $\chi : L^2(\Omega)\to [0,r]$. Here $\Theta^\ell
(\varphi)\equiv \chi (\varphi(-r))$   and $ \Theta^u (\varphi)=r$.
It is easy to see that the above delay function $\eta$
(\ref{sdd12-4-4}) ignores values of $\varphi$ at points $\theta\in
(-r,-\chi(\varphi(-r)))\cup (-\chi(\varphi(-r)),0]$ and uses just
two values of $\varphi$ at points $\theta=-r, \,
\theta=-\chi(\varphi(-r))$. In our notations, the delayed segment
$[-\Theta^u  (\varphi), -\Theta^\ell (\varphi)]=
[-r,-\chi(\varphi(-r)]$ is state-dependent.

In the same way, one has
$$\eta(\varphi)= \sum^N_{k=1} p_k\left( \varphi(-\chi^k(\varphi(-r))\right) \hbox{ with } p_k,\chi^k: L^2(\Omega)\to
[0,r].
%\quad \min r_k >0.
$$
In this case $[-\Theta^u  (\varphi), -\Theta^\ell (\varphi)]=
[-r,-\min_k\{\chi^k(\varphi(-r))\}]$. A slightly more general
example is
$$\eta(\varphi)= \sum^N_{k=1} p_k\left( \varphi(-\chi^k(\varphi(-r^k))\right) \hbox{ with } p_k,\chi^k: L^2(\Omega)\to
[0,r], \quad \min r^k \in (0,r].
$$
Here $\Theta^u
(\varphi)=\max\{r^1,\ldots,r^N,\chi^1(\varphi(-r^1)),\ldots,\chi^N(\varphi(-r^N))\}
$ and

$\Theta^\ell
(\varphi)=\min\{r^1,\ldots,r^N,\chi^1(\varphi(-r^1)),\ldots,\chi^N(\varphi(-r^N))\}$.

Examples of integral delay terms are as follows
$$\eta(\varphi)= \int^{-\chi^1(\varphi(-r^1))}_{-\chi^2(\varphi(-r^2))} p_1(\varphi(\theta)) g(\theta)\, d\theta,\quad
\hbox{ and } \quad\eta(\varphi)= p_1\left(
\int^{-\chi^1(\varphi(-r^1))}_{-\chi^2(\varphi(-r^2))}
\varphi(\theta)  g(\theta)\,
 d\theta\right).%\quad \eta_{ign}>0,  \hbox{ etc. }
 $$
 Similar to the previous example, $\Theta^u
(\varphi)=\max\left\{r^1,r^2,\chi^1(\varphi(-r^1)),\chi^2(\varphi(-r^2))\right\}
$ and

$\Theta^\ell
(\varphi)=\min\left\{r^1,r^2,\chi^1(\varphi(-r^1)),\chi^2(\varphi(-r^2))\right\}$.
%To satisfy assumption (H1) for the functions given above it is
%sufficient to assume that $\inf p_i(\cdot) >0. $
 }%

\medskip

{\bf Remark~8.} {\it It is interesting to notice that an assumption
similar  to the existence of upper function $\Theta^u (\cdot)$ is
used in \cite{Walther-preprint-2010} for ODEs with SDD (locally
bounded delay). On the other hand, an assumption similar to (H) is
used in \cite{Hartung-preprint-2010} for neutral ODEs (see (A4)(ii)
in \cite{Hartung-preprint-2010}), but together with another
assumption on SDD to be bounded from below by a constant $r_0 > 0$
(c.f. remark~6).

}%

\medskip

Following \cite[theorem~1]{Rezounenko_NA-2009} we have the first
result

\medskip

\medskip

{\bf Theorem~1.} {\it Let both upper and low %ignoring
functions $\Theta^u , \Theta^\ell : C\to [0,r]$ be continuous and
$\Theta^\ell (\varphi)>0, \forall \varphi\in C$. Assume the delay
function $\eta : C\to [0,r]\subset R_{+}$ is continuous and
satisfies assumption $(\widehat{H})$;
%the function $b : R\to R$ is a %%%%%locally
%Lipschitz map, satisfying $|b(w)|\le C_1 |w|+C_b$ with $C_i\ge 0$; % the assumptions (H1) and (H2)
% $f: \Omega -\Omega \to R$ is a bounded and measurable function ($|f(z )|\le M_f,
%\forall z\in \Omega -\Omega$).
the mapping $B$ is Lipschitz continuous (see (\ref{sdd12-2-2})).

Then for any
  initial function $\varphi\in C,$ initial value problem (\ref{sdd8-1}), (\ref{sdd8-ic})
  has an {\tt unique} mild solution $u(t), t\ge 0$ (given by
  proposition~1).

  If we define the {\tt evolution operator} $S_t: C\to C$ by the
formula $S_t \varphi \equiv u_t,$ where $u(t)$ is the unique mild
solution of (\ref{sdd8-1}), (\ref{sdd8-ic}) with initial function
$\varphi$, then the pair $(S_t, C)$ constitutes a dynamical system
i.e. the following properties are satisfied:
\begin{enumerate}
 \item $S_0=Id$ ( identity operator in $C$ );
 \item $\forall\,\, t,\tau \ge 0\quad  \Longrightarrow \quad  S_t\, S_\tau = S_{t+\tau}$;
 \item $t\mapsto S_t$ is a strongly continuous in $C$ mapping;
 \item for any $t\ge 0$ the evolution operator $S_t$ is continuous in $C$ i.e. for any
 $\{\varphi^n\}^\infty_{n=1}\subset C$ such that $||\varphi^n -\varphi||_C\to 0$ as
 $n\to \infty$, one has $||S_t\varphi^n -S_t\varphi||_C\to 0$ as
 $n\to \infty.$
\end{enumerate}
}
\medskip

The proof follows the line of \cite[theorem~1]{Rezounenko_NA-2009}
taking into account that condition $\Theta^\ell (\varphi)>0, \forall
\varphi\in C$ implies that for any fixed $\varphi\in C$, due to the
continuity of $\Theta^\ell :C\to [0,r]$, there exists a neighbourhood %neighborhood
$U(\varphi)\subset C$ such that for all $\psi\in U(\varphi)$ one has
$\Theta^\ell (\psi)\ge {1\over 2}\Theta^\ell (\varphi)> 0.$ That
means that in $U(\varphi)\subset C$ we have the (state-independent)
%ignoring
 condition (H) with $\eta_{ign}={1\over 2}\Theta^\ell
(\varphi)> 0 $ and all the arguments presented in
\cite[theorem~1]{Rezounenko_NA-2009} could be directly applied to
this case. \hfill \finproof

\medskip

{\bf Remark~9.} {\it We do {\tt not} assume that the upper and low
%ignoring
 functions $\Theta^u , \Theta^\ell $ (which are used in
$(\widehat{H})$ to present the delayed segment $[-\Theta^u
(\varphi), -\Theta^\ell (\varphi)]$) are the functions presenting
the smallest possible delayed segment. More precisely, it is
possible that there exist two other functions $\widetilde{\Theta}^u,
\widetilde{\Theta}^\ell$ such that for all $\varphi\in C$ one has
$0\le \Theta^\ell (\varphi) \le \widetilde{\Theta}^\ell(\varphi) \le
\widetilde{\Theta}^u(\varphi) \le \Theta^u  (\varphi)\le r $ and the
same delay    $\eta$ satisfies $(\widehat{H})$ with
$\widetilde{\Theta}^u, \widetilde{\Theta}^\ell$ as well.
}%

\medskip

Our next step in studying the {\it state-dependent} %ignoring
condition $(\widehat{H})$ is an attempt to avoid the condition
$\Theta^\ell (\varphi)>0, \forall \varphi\in C$. We are going to
consider the general case $\Theta^\ell (\varphi)\ge 0, \forall
\varphi\in C$ with a {\tt non-empty} set $Z\equiv \{ \varphi\in C :
\Theta^\ell (\varphi) = 0\} \neq \emptyset$.

\medskip

{\bf Theorem~2.} {\it Assume
%the function $b : R\to R$ is a %%%%%locally
%Lipschitz map, satisfying $|b(w)|\le C_1 |w|+C_b$ with $C_i\ge 0$; % the assumptions (H1) and (H2)
% $f: \Omega -\Omega
%\to R$ is a bounded and measurable function ($|f(z )|\le M_f,
%\forall z\in \Omega -\Omega$).
the mapping $B$ is Lipschitz continuous (see (\ref{sdd12-2-2})).

\medskip

Moreover, let the following conditions be satisfied:

1) both upper and low %ignoring
functions $\Theta^u , \Theta^\ell :
C\to [0,r]$ are continuous; %%%%$\Theta^\ell (\varphi)\ge 0,\forall \varphi\in C$; %
\par 2)    $Z\equiv \{ \varphi\in C :
\Theta^\ell (\varphi) = 0\} \subset C{\cal L}_L\equiv \left\{  \varphi\in C : \sup\limits_{t\neq s} \frac{||\varphi(t)-\varphi(s) ||}{|t-s|} \le L\right\} $; %
\par 3)  delay function $\eta : C\to [0,r]\subset R_{+}$ is
continuous and satisfies assumption $(\widehat{H})$;
\par 4) $\forall \varphi\in Z \Rightarrow  \eta (\varphi) > 0$;
\par 5) $\exists U_\omega (Z)\equiv \{ \chi\in C : \exists\nu\in Z : ||\chi-\nu||_C\le \omega\}, \exists L_\eta >0 :
\forall \varphi,\psi\in U_\omega (Z)   \Rightarrow  $
%\par $
$$|\eta (\varphi) -\eta(\psi)|\le L_\eta\cdot ||\varphi-\psi||_C.$$

\medskip

Then for any
  initial function $\varphi\in C,$ initial value problem (\ref{sdd8-1}), (\ref{sdd8-ic})
  has an {\tt unique} mild solution $u(t), t\ge 0$ (given by
  proposition~1). Moreover, the pair $(S_t, C)$ constitutes a dynamical system (see thm~1).
}%

\medskip

{\it Proof of theorem~2.} Let us consider $\varphi\in C$ which is an
initial condition (see (\ref{sdd8-ic})). We start with the simple
case $\varphi\not\in Z.$ By definition of Z, we have $\Theta^\ell
(\varphi)>0. $ Hence we apply the same arguments as in the proof of
theorem~1 (the state-independent %ignoring
condition (H) is satisfied
locally).

The rest of the proof is devoted to the case $\varphi\in Z.$ We
remind some estimates similar to estimates (6)-(13)  in
\cite{Rezounenko_NA-2009}. Denote by $u^k(t)$ any solution of
(\ref{sdd8-1}),(\ref{sdd8-ic}) with the initial function $\varphi^k$
and by $u(t)$ any solution of (\ref{sdd8-1}),(\ref{sdd8-ic}) with
the initial function $\varphi.$

We use the variation of constants formula for parabolic equation
%(\ref{sdd8-4-01})
(with $\widetilde A \equiv A+d\cdot E$)

\begin{equation}\label{sdd8-4-1}
u(t)=e^{-\widetilde A t}u(0) + \int^{t}_0 e^{-\widetilde A (t-\tau)}
%\int_\Omega b(u(\tau-\eta (u_\tau), y)) f(\cdot -y) dy
B(u(\tau- \eta(u_\tau))) \, d\tau,
\end{equation}
\begin{equation}\label{sdd8-4-2}
u^k(t)=e^{-\widetilde A t}u^k((0) + \int^{t}_0 e^{-\widetilde A
(t-\tau)}
%\int_\Omega b(u^k(\tau-\eta (u^k_\tau), y)) f(\cdot -y) dy
B(u^k(\tau- \eta(u^k_\tau))) \, d\tau.
 \end{equation}
Using $||e^{-\widetilde A t}||\le 1$ and  $||e^{-\widetilde A
(t-\tau)}||\le 1$, we get

$$||u^k(t)-u(t) ||\le ||u^k(0)-u(0)|| %$$ $$
+ \int^{t}_0
||
%\int_\Omega \left[ b(u^k(\tau-\eta (u^k_\tau), y)) - b(u(\tau-\eta
%(u_\tau), y)) \right] f(\cdot -y) dy
B(u^k(\tau- \eta(u^k_\tau)))- B(u(\tau- \eta(u_\tau))) %
 || \, d\tau $$
\begin{equation}\label{sdd8-4-3}
= ||\varphi^k(0)-\varphi(0)|| + { J}^k_1(t) + { J}^k_2(t),
\end{equation}
where we denote (for $s\ge 0, x\in \Omega$)
\begin{equation}\label{sdd8-4-4}
{ J}^k_1(s)\equiv { J}^k_1(s)(x)\equiv \int^{s}_0 ||
%\int_\Omega
%\left[ b(u^k(\tau-\eta (u^k_\tau), y)) - b(u(\tau-\eta (u^k_\tau),
%y)) \right] f(x -y) dy
B(u^k(\tau- \eta(u^k_\tau)))- B(u(\tau- \eta(u^k_\tau)))||  \,
d\tau,
\end{equation}
\begin{equation}\label{sdd8-4-5}
{ J}^k_2(s)\equiv { J}^k_2(s)(x)\equiv\int^{s}_0 ||
%\int_\Omega
%\left[ b(u(\tau-\eta (u^k_\tau), y)) - b(u(\tau-\eta (u_\tau), y))
%\right] f(x -y) dy
B(u(\tau- \eta(u^k_\tau)))- B(u(\tau- \eta(u_\tau))) %
|| \, d\tau.
\end{equation}
Using the Lipschitz property (\ref{sdd12-2-2}) of $B$, %$b$,
one easily gets %\marginpar{???}
$${ J}^k_1(t)\le L_B %M_f |\Omega| L_b
\int^{t}_0 ||u^k(\tau- \eta(u^k_\tau))- u(\tau- \eta(u^k_\tau))|| \,
d\tau
$$
\begin{equation}\label{sdd8-4-6}  \le L_B %M_f |\Omega| L_b
t \max_{s\in [-r,t]} ||u^k(s)- u(s)||.
\end{equation}
Estimates (\ref{sdd8-4-6}), (\ref{sdd8-4-3}) and property ${
J}^k_2(s) \le { J}^k_2(t)$ for $s\le t \le t_0$ %<\eta_{ign}$
 give
$$\max_{t\in [0,t_0]} ||u^k(t)- u(t)|| \le
||\varphi^k(0)-\varphi(0)|| + L_B %M_f |\Omega| L_b %
t_0 \max_{s\in [-r,t_0]} ||u^k(s)- u(s)|| + { J}^k_2(t_0).
$$ Hence
\begin{equation}\label{sdd8-4-6a}
\max_{s\in [-r,t_0]} ||u^k(s)- u(s)|| \le ||\varphi^k-\varphi||_C +
L_B %M_f |\Omega| L_b
t_0 \max_{s\in [-r,t_0]} ||u^k(s)- u(s)|| + {
J}^k_2(t_0).
\end{equation}

%We choose $t_0< [M_f |\Omega| L_b]^{-1}$ (to satisfy $M_f |\Omega|
%L_b t_0 <1$) and get
%\begin{equation}\label{sdd8-4-6a}
%(1-M_f |\Omega| L_b t_0) \max_{s\in [-r,t_0]} ||u^k(s)- u(s)|| \le
%||\varphi^k-\varphi||_C  +  { J}^k_2(t_0).
%\end{equation}

%Our goal is to show that ${ J}^k_2(t_0) \to 0$ as $k\to \infty.$
Now we study properties of ${ J}^k_2$ which essentially {\tt differ}
from the ones in \cite{Rezounenko_NA-2009} since (H) is not
satisfied. The Lipschitz property of $B$ %$b$
implies
\begin{equation}\label{sdd8-4-7}
{ J}^k_2(t_0)\le L_B %M_f |\Omega|  L_b
\int^{t_0}_0 || u(\tau-\eta
(u^k_\tau)) - u(\tau-\eta (u_\tau)) || \, d\tau.
\end{equation}

Since $\varphi\in Z$, property 4) gives $\eta(\varphi)>0$. Due to
the continuity of $\eta$ (see 3)), %there exists a neighborhood %neighbourhood
\begin{equation}\label{sdd12-4-1}
\exists U_\alpha(\varphi)\equiv \left\{ \psi\in C : ||\varphi-\psi
||_C\le \alpha\right\} : \forall\psi \in
U_\alpha(\varphi)\Rightarrow \eta(\psi) \ge {3\over 4}
\eta(\varphi)>0.
\end{equation}
We choose $\alpha<\omega$ (see property 5). %
 By definition, a solution is strongly continuous function
(with values in $L^2(\Omega)$), hence for any two solutions $u(t)$
and $u^k(t)$ there exist two time moments $t_\varphi,
t_{\varphi^k}>0$
%\in (0,t_0]$
such that for all $t\in (0,t_\varphi]$ one has $u_t\in
U_\alpha(\varphi)$ and for all $t\in (0,t_{\varphi^k}]$ one has
$u^k_t\in U_\alpha(\varphi)$.

\medskip

{\bf Remark~10.} {\it More precisely, we assume that $\exists
N_\alpha \in N$ such that for all $k\ge N_\alpha$ one has $\varphi^k
\in U_{\alpha/2}(\varphi)$ and hence there exists time moment
$t_{\varphi^k} \in (0,t_0]$ such that for all $t\in
(0,t_{\varphi^k}]$ one has $u^k_t\in U_\alpha(\varphi)$. The last
assumption ($ \exists N_\alpha \in N : \forall k\ge N_\alpha
\Rightarrow \varphi^k \in U_{\alpha/2}(\varphi)$) is not restrictive
since for the uniqueness of solutions we have $\varphi^k =\varphi$
while for the continuity with respect to initial data (see below) we
have $\varphi^k \to \varphi$ in $C$.
}%

\medskip

{\bf Remark~11.} {\it It is important to notice that we take {\tt
any} solution from the set of solutions of IVP
(\ref{sdd8-1}),(\ref{sdd8-ic}) with the initial function $\varphi$
(and denote it by $u(t)$) and take {\tt any} solution from the set
of solutions of IVP (\ref{sdd8-1}),(\ref{sdd8-ic}) with the initial
function $\varphi^k$ (and denote it by $u^k(t)$) i.e. the values
$t_\varphi, t_{\varphi^k}$ may depend on the choice of these {\tt
two} solutions.}

\medskip

These and (\ref{sdd12-4-1}) imply that  for all $\tau\in [0, t_1]$,
with $ t_1 \le \min\{ t_\varphi; t_{\varphi^k}; {3\over 4}
\eta(\varphi) \} $ %
%$ t_1\equiv\min\{ t_\varphi; t_{\varphi^k}; t_0;{3\over 4}
%\eta(\varphi) \} $ %
 one gets $\tau-\eta (u_\tau)\le 0, \tau-\eta
(u^k_\tau) \le 0$ and $u(\tau-\eta (u_\tau))= \varphi(\tau-\eta
(u_\tau)), u(\tau-\eta (u^k_\tau)) = \varphi(\tau-\eta (u^k_\tau))
$. Hence, see (\ref{sdd8-4-7}) and properties 2), 5),
$${ J}^k_2(t_1)\le L_B %M_f |\Omega|  L_b
\int^{t_1}_0 || \varphi(\tau-\eta
(u^k_\tau)) - \varphi(\tau-\eta (u_\tau)) || \, d\tau \le
L_B %M_f|\Omega|  L_b
L \int^{t_1}_0 |\eta (u^k_\tau) - \eta (u_\tau) | \,
d\tau
$$
$$ \le L_B %M_f |\Omega|  L_b
L L_\eta t_1  \max_{s\in [-r,t_1]} ||u^k(s)- u(s)||.
$$
Finally, we get (see the last estimate and (\ref{sdd8-4-6a}))
$$(1-L_B %M_f |\Omega| L_b
t_1 [1+L
L_\eta ]) \max_{s\in [-r,t_1]} ||u^k(s)- u(s)|| \le
||\varphi^k-\varphi||_C. $$

Choosing small enough $t_1>0$ (to have $1-L_B %M_f |\Omega| L_b
t_1 [1+L
L_\eta ]>0$) i.e.
\begin{equation}\label{sdd12-4-2}
t_1\equiv\min\left\{ t_\varphi; t_{\varphi^k}; %t_0;
{3\over 4} \eta(\varphi);q L_B %(M_f |\Omega| L_b %
 [1+L L_\eta ])^{-1}
\right\} \quad \hbox{ for any fixed }\quad  q\in (0,1),
\end{equation}%
 we get
\begin{equation}\label{sdd12-4-3}
\max_{s\in [-r,t_1]} ||u^k(s)- u(s)|| \le (1-L_B %M_f |\Omega| L_b
t_1 [1+L L_\eta ])^{-1} ||\varphi^k-\varphi||_C.
\end{equation}
 It is easy to see that (\ref{sdd12-4-3}) particularly implies the {\it uniqueness} of mild solutions to
I.V.P. (\ref{sdd8-1}),(\ref{sdd8-ic}) in case when
$\varphi^k=\varphi.$

\medskip

It gives us the possibility to define the {\tt evolution operator}
$S_t: C\to C$ by the formula $S_t \varphi \equiv u_t,$ where $u(t)$
is the unique mild solution of (\ref{sdd8-1}), (\ref{sdd8-ic}) with
initial function~$\varphi$.

\medskip

Our next goal is to prove that pair $(S_t, C)$ constitutes a {\tt
dynamical system} (see the properties $1.-4.$ as they are formulated
 in theorem~1). As in \cite[p.3981]{Rezounenko_NA-2009}, properties
1, 2 are consequences of the uniqueness of mild solutions. Property
3 is given by Proposition~1 since the solution is a continuous
function $u\in C([-r,T]; L^2(\Omega)).$

Let us prove property 4. We consider any sequence $\{
\varphi^k\}^\infty_{k=1}\subset C$, which converges (in space $C$)
to $\varphi.$ Denote by $u^k(t)$ the (unique!) mild solution of
(\ref{sdd8-1}),(\ref{sdd8-ic}) with the initial function $\varphi^k$
and by $u(t)$ the (unique!) mild  solution of
(\ref{sdd8-1}),(\ref{sdd8-ic}) with the initial function $\varphi.$

One could think that (\ref{sdd12-4-3}) already provides the
continuity with respect to initial data, but there is an important
technical property used in developing (\ref{sdd12-4-3}) i.e. the
choice of $t_1$ (see (\ref{sdd12-4-2}) and remark~11). In contrast
to the previous study, now we have infinite set of functions $\{
\varphi^k\}^\infty_{k=1}\subset C$, so it may happen that $t_1=t^k_1
\to 0$ when $k\to \infty.$

We remind (see the text after (\ref{sdd12-4-1})) that {\tt two} time
moments $t_\varphi, t_{\varphi^k}>0$ % \in (0,t_0]$
have been chosen such that for all $t\in (0,t_\varphi]$ one has
$u_t\in U_\alpha(\varphi)$ and for all $t\in (0,t_{\varphi^k}]$ one
has $u^k_t\in U_\alpha(\varphi)$. Now our goal is to show that
infinite number of moments $t_\varphi$,
$\{t_{\varphi^k}\}^\infty_{k=1}$ could be chosen in such a way that
$t_2\equiv \inf_{k\in N} \{t_\varphi,
t_{\varphi^k} \} %^\infty_{k=1}
>0$ and $u_t, u^k_t \in U_\alpha(\varphi) $ for all $t\in (0,t_2]$.
To get this, we use the standard proof of the existence of a mild
solution by a fixed point argument (see e.g. \cite[p.46, thm
2.1]{Wu_book}). More precisely, let $U$ be an open subset of $C$ and
$\widetilde{F} : [0,b]\times U \to L^2(\Omega)$ be continuous. For
$\varphi\in C$ and any $y\in Y_1\equiv \{ y\in
C([-r,t_3];L^2(\Omega)) : y(0)=\varphi(0)) \}$ we consider the
extension function $\hat y$ as follows
$$
\hat y(s)\equiv \left[\begin{array}{ll}
  \varphi(s) & \hbox{ for } s \in [-r, 0]; \\
  y(t) & \hbox{ for } s\in (0, t_3] \\
\end{array}.
\right. $$ Let $Y_2\equiv \{ y\in Y_1 : \hat y_t\in
\overline{B_\delta(\varphi)}\hbox{ for } t\in [0,t_3]\}.$ Consider a
mapping $G$ on $Y_2$ as follows
$$G(y)(t)\equiv e^{-\widetilde A t}\varphi(0) + \int^{t}_0 e^{-\widetilde A (t-\tau)}
\widetilde{F}(\hat y_\tau) \, d\tau.
$$
One can check (see \cite[p.46,47, thm 2.1]{Wu_book}), that $G$ maps
$Y_2$ into $Y_2$ provided $t_3\equiv \min \{ t^\prime; b;
\delta/(3N);\delta\}$. Here we use notations of \cite[p.46]{Wu_book}
chosen as follows. Constants $\delta>0$ and $N>0$ are such that
$||\widetilde{F}(\psi)||\le N$ for all $\psi\in
\overline{B_\delta(\varphi)}\equiv \{ \psi\in C :
||\psi-\varphi||_C\le \delta\}$, $||e^{-\widetilde A t}||\le M=1$.
The time moment $t^\prime<r$ is chosen so that if $0\le t\le
t^\prime$ then $||\varphi(t+\theta)-\varphi(\theta)||< \delta/3$ and
$||e^{-\widetilde A t}\varphi(0) - \varphi(0)|| <\delta/3$. The
solution is given by a fixed point $y=G(y)$. For our goal it is
sufficient to choose $\delta\le\alpha$ and $t_2\le t_3$ to get $u_t,
u^k_t \in U_\alpha(\varphi) $ for all $t\in (0,t_2]$. Here we use
$\varphi^k$ instead of $\varphi$ when necessary.  The crucial point
here is the possibility to choose $t^\prime$ (and hence $t_3$ and
$t_2$) independent of $k\in N$. The choice of $t^\prime<r$  so that
if $0\le t\le t^\prime$ then $||\varphi(t+\theta)-\varphi(\theta)||<
\delta/3$ and $||\varphi^k(t+\theta)-\varphi^k(\theta)||< \delta/3$
for {\tt all} $k\in N$ is possible due to the convergence of
$\varphi^k$ (to $\varphi$ in $C$). Since any convergent sequence is
a pre-compact set in $C$, the desired property is the equicontinuity
given by the Arzela-Ascoli theorem. Now estimate (\ref{sdd12-4-3})
can be applied to our case and this completes the proof of property
4 and theorem~2. \finproof

 \medskip

Discussing assumptions of theorem~2, let us present a constructive
{\it example} of low %ignoring
function $\Theta^\ell $ which satisfies assumption 2). Consider any
compact and convex set $K_C\subset C{\cal L}_L\subset C$. For
example, for any compact and convex set $K\in L^2(\Omega)$, the set
$\{ \varphi\in C : \varphi\in C{\cal L}_L, \forall \theta\in [-r,0]
\Rightarrow \varphi(\theta)\in K\}$ is compact (by Arzela-Ascoli
theorem) and convex. First, constructing $\Theta^\ell $, we set
$\Theta^\ell (\varphi)=0$ for all $\varphi\in K_C$. Second, we take
any $p\in (0,r]$ and set $\Theta^\ell (\varphi)=p$ for all
$\varphi\in C$ such that $\hbox{dist}_C (\varphi,K_C) \ge 1$. Third,
for any $\varphi\in C $ such that $\hbox{dist}_C (\varphi,K_C) \in
(0,1)$ we find an unique $\widehat \varphi\in K_C$ such that
$\hbox{dist}_C (\varphi,K_C) =||\varphi-\widehat \varphi||_C$. Such
$\widehat \varphi\in K_C$ exists by the classical Weierstrass
theorem since $f(\psi)\equiv \hbox{dist}_C (\varphi,\psi) : K_c \to
(0,1)$ is continuous ($\varphi$ is fixed) and $K_C$ is compact. The
uniqueness of $\widehat \varphi$ follows from the convexity of
$K_C$. Finally, we set $\Theta^\ell (\varphi)= p\cdot \hbox{dist}_C
(\varphi,\widehat \varphi)\in (0,p)$ for all $\varphi\in C :
\hbox{dist}_C (\varphi,K_C) \in (0,1)$. By construction,
$\Theta^\ell $ satisfies 2).

\medskip

As for asymptotic behavior, we study of the long-time %%%%%asymptotic
behavior of the dynamical system $( S_t,C )$, constructed in
theorems~1 and 2. Similar to \cite[theorem~2]{Rezounenko_NA-2009} we
have the following result.

\medskip

{\bf Theorem~3.} {\it Assume all the assumptions of theorems~1 or 2
are satisfied and additionally mapping $B$ (see (\ref{sdd12-2-1}))
is bounded.
 Then the dynamical
system $( S_t,C )$ has a compact global attractor ${\cal A}$ which
is a compact set in all spaces $C_\delta\equiv C([-r,0];
D(A^\delta)), \forall\delta\in [0,{1\over 2}).$}

\medskip

{\bf Lemma.} {\it Let all the assumptions of theorem~2 be satisfied.
Then the global attractor ${\cal A}$ (see theorem~3) is a subset of
$C{\cal L}_{\widetilde{L}}$ (c.f. condition 2 in theorem~2).
}%

\medskip

{\bf Remark~12. } {\it Lemma gives a possibility to consider system
(\ref{sdd8-1}), (\ref{sdd8-ic}) with a state-dependent delay
function $\eta$ which does not ignore values of its argument
$\varphi$ for all points $\varphi\in {\cal A}$ i.e. no information
is lost %(ignored )
 on the  global attractor ${\cal A}$.
}%

\medskip

{ } {\it Proof of lemma.} Consider any solution $u_t\in {\cal A}$.
 Let us denote $f(t)\equiv
F(u_t)$ and prove that $f$ is H\"older continuous.

We will need the following property, proved in \cite[estimate (29)
with $\delta=0$]{Rezounenko_NA-2009}
\begin{equation}\label{sdd12-4-5}
||u(t_1)-u(t_2)||\le  L_0 |t_1-t_2|^{1/2}
\end{equation}
for any solution, belonging to the ball of dissipation
(particularly, for any solution belonging to the attractor). Here
$L_0$ is independent of solution $u$.

One can check that
\begin{equation}\label{sdd12-4-6}
||f(t_1)-f(t_2)||\le L_B %M_f |\Omega| L_b
\cdot ||u(t_1-\eta (u_{t_1})) -
u(t_2-\eta (u_{t_2}))||.
\end{equation}

Using (\ref{sdd12-4-5}), the Lipschitz property of $\eta$ (see 5 in
theorem~2), we get from (\ref{sdd12-4-6}) that
$$||f(t_1)-f(t_2)||\le L_B %M_f |\Omega| L_b
L_0 \cdot |t_1-\eta (u_{t_1})-(t_2-\eta (u_{t_2}))|^{1/2}
$$
$$\le L_B %M_f |\Omega| L_b
L_0 \cdot \left(|t_1-t_2| + |\eta (u_{t_1})-\eta (u_{t_2}) |\right)^{1/2}
\le \hbox{[ using 5 in theorem~2 ]} \le
$$
$$\le L_B %M_f |\Omega| L_b
L_0 \cdot \left(|t_1-t_2| + L_\eta ||u_{t_1}-u_{t_2}|| \right)^{1/2}
\le \hbox{[ using (\ref{sdd12-4-5}) ]} \le
$$
$$\le L_B %M_f |\Omega| L_b
L_0 \cdot \left(|t_1-t_2| + L_\eta L_0 |t_1-t_2|^{1/2}
\right)^{1/2}.
$$
$$\le L_B %M_f |\Omega| L_b
L_0 \cdot \left(|t_1-t_2|^{1/2} + (L_\eta L_0)^{1/2} |t_1-t_2|^{1/4}
\right).
$$

Finally, for $|t_1-t_2|<1$ one has
\begin{equation}\label{sdd12-4-7}
||f(t_1)-f(t_2)||\le L_B %M_f |\Omega| L_b
L_0 \cdot \left\{ 1+(L_\eta L_0)^{1/2}\right\} |t_1-t_2|^{1/4}.
\end{equation}

Let us consider $\forall \psi\in {\cal A}$. It is well-known that
the attractor consists of whole trajectories i.e. $u_s\in {\cal A},
\, \forall s\in R$. We take any $t_0>r>0$ and get $\varphi\in {\cal
A}$ such that $S_{t_0}\varphi =\psi.$ Consider the variation of
constants formula for parabolic equation (with $\widetilde A \equiv
A+d\cdot E$ see (\ref{sdd8-4-1}))
\begin{equation}\label{sdd12-4-8}
u(t)=e^{-\widetilde A t}\varphi(0) + \int^{t}_0 e^{-\widetilde A
(t-\tau)} F(u_\tau) \, d\tau.
\end{equation}

The first term in the above formula (\ref{sdd12-4-8}) is Lipschitz
for $t>t_0$ due to the standard estimate $||e^{-\widetilde At_1}v -
e^{-\widetilde At_2}v || \le (t_1 e)^{-1} ||v|| \cdot |t_1-t_2|, \,
0<t_1<t_2.$ Moreover it is uniformly Lipschitz for any
$v=\varphi(0), \varphi\in {\cal A}$ since $||e^{-\widetilde
At_1}\varphi(0), -  e^{-\widetilde At_2}\varphi(0), || \le (r
e)^{-1} ||v|| \cdot |t_1-t_2| \le (r e)^{-1} C(0) \cdot |t_1-t_2|,
\, r<t_1<t_2.$ Here $||\varphi(0)||\le C(0)$ due to the
dissipativeness of the dynamical system $(S_t;C)$ (for more details
see \cite[estimate (23)]{Rezounenko_NA-2009}).

To prove that the second term in (\ref{sdd12-4-8}) is Lipschitz for
$t>t_0$ we need the following

\medskip

{\bf Proposition } \cite[lemma 3.2.1]{Henry-1981}. {\it Let
$\widetilde A$ be a sectorial operator in Banach space $X$. Assume
function $f:(0,T)\to X$ is locally H\"{o}lder continuous and
$\int^\rho_0 ||f(s)||_X\, ds <\infty$ for some $\rho>0.$ Denote by
$\Phi(t)\equiv \int^t_0 e^{-\widetilde A(t-s)}f(s)\, ds.$ Then
function $\Phi(\cdot)$ is continuous on $[0,T)$, continuously
differentiable on $(0,T)$, $\Phi(t)\in D(\widetilde A)$ for $0<t<T$
and $d\Phi(t)/dt + \widetilde A \Phi(t)=f(t)$ for $0<t<T$ and
$\Phi(t)\to 0$ in $X$ as $t\to 0+$.
}%

\medskip

{\bf Remark~13. } {\it Our operator $\widetilde A$ is sectorial
since any self-adjoint densely defined bounded from below operator
in a Hilbert space is sectorial (see e.g. \cite[example~2,
p.26]{Henry-1981}).
}%

\medskip

We apply the above proposition to $f(t)\equiv F(u_t)$ and use
(\ref{sdd12-4-7}). The property $\int^\rho_0 ||f(s)||_X\, ds
<\infty$ for some $\rho>0$ follows from the dissipativeness
$||u(t)||\le C(0)$, the continuity of $F:C\to L^2(\Omega)$ and the
strong continuity of mild solution $u$. One uses the continuous
differentiability of $\Phi$ on $[t_0-r,t_0]\subset (0,T)$ which
implies that $\max_{t\in [t_0-r,t_0]} ||\Phi^\prime (t)|| \equiv
M_{\Phi;1}<\infty.$ In our case $\Phi$ represents the second term in
(\ref{sdd12-4-8}) which is proved to be Lipschitz continuous with
Lipschitz constant $M_{\Phi;1}$ independent of $u$. The proof of
lemma is complete.\hfill \finproof

\medskip%%%%%%%%%%%%%%%%%%%%%%%%%%%%%%%%%%%%%%%%%%%%%%%%%%%%%%%%%%%%%%%%%%%%%%%%%%%%%%%%%%%%%%%

{\bf Remark~14. } {\it %All the results above are valid for
One can also easily extend the method developed here to the case of
non-autonomous nonlinear delay terms, for example, using another
nonlinear function $\hat b :R\times R \to R$ (see remark 2) instead
of $b$ to have $\big( \hat F(t,u_t) \big)(x)=\hat b(t,u(t-\eta
(u_t),x))$ or $\big( \hat F(t,u_t) \big)(x)=\int_\Omega \hat
b(t,u(t-\eta (u_t), y)) f(x-y) dy $ in equation (\ref{sdd8-1}).
}%

%%%%%%%%%%%%%%% %%%%%%%%%%%%%%%%%% %%%%%%%%%%%%%%%% %%%%%%%%%%%%%
\medskip

As an application we can consider the diffusive Nicholson's
blowflies equation (see e.g. \cite{So-Yang} %,So-Wu-Yang})
   with
state-dependent delays. More precisely, we consider equation
(\ref{sdd8-1}) where $-A$ is the Laplace operator with the Dirichlet
boundary conditions, $\Omega\subset R^{n_0}$ is a bounded domain
with a smooth boundary, the function $f$ (see remark~2) can be,
 for example, $ f(s)={1\over \sqrt{4\pi\alpha}}
e^{-s^2/4\alpha}$, as in \cite{So-Wu-Zou} (for the non-local in
space variable nonlinearity) or Dirac delta-function to get the
local in space variable nonlinearity, the nonlinear function $b$ is
given by $b(w)=p\cdot we^{-w}.$ Function $b$ is bounded%and $b(w)>0$ for all $w\neq 0.$
, so for any continuous delay function $\eta$, satisfying ($\widehat
H$), the conditions of theorems~1,2 are valid. As a result, we
conclude that the initial value problem
(\ref{sdd8-1}),(\ref{sdd8-ic}) is well-posed in $C$ and the
dynamical system $(S_t,C)$ has a global attractor (theorem~3).
\medskip

%\noindent {\bf Acknowledgements.}  The author wishes to thank
%Professors I.D.~Chueshov and  %\\
%H.-O.~Walther for useful discussions of an early version of the
%manuscript.
%%% He is also thankful to anonymous referee for very useful
%%% comments which lead to a better presentation of the results.

\medskip

{\bf Acknowledgement}. The author wishes to thank I.D. Chueshov and
H.-O.Walther for useful discussions of an early version of the
manuscript.

\medskip

\hfill Kharkiv, November 19, 2010

\medskip

\hfill {\it E-mail}: rezounenko@univer.kharkov.ua

%
% $\eta^\ell, \quad  \rho^\ell, \quad  \Theta^\ell, \quad [-\Theta^u(\varphi), -\Theta^\ell(\varphi)],
%\quad \Theta^u \quad\eta^u, \quad \rho^u, \quad$
%

\end{document}